\documentclass[graybox]{svmult}

\usepackage{mathptmx}       
\usepackage{helvet}         
\usepackage{courier}        
\usepackage{type1cm}        
\usepackage{makeidx}         
\usepackage{graphicx}        
\usepackage{multicol}        
\usepackage[bottom]{footmisc}

\makeindex             

\usepackage{amsmath}
\usepackage{amssymb}
\usepackage{amsfonts}
\usepackage[]{algorithm}
\usepackage[noend]{algpseudocode}
\usepackage{color}
\usepackage{array}
\usepackage{psfrag}
\usepackage{epsfig}
\usepackage{subfigure}
\usepackage{float}
\usepackage{theorem}
\usepackage{multirow}
\usepackage{enumitem}  

\def\div{\nabla\cdot}

\def\n{{\bm n}}

\def\to{t^{n+1}}

\def\E{\mathbf{E}}

\def\R{\hbox{\doble R}}

\def\cl {\nonumber \\}
\def\el {\nonumber }

\font\doble=msbm10 scaled\magstep1             
\newcommand{\beqn}{\begin{eqnarray}}
\newcommand{\eeqn}{\end{eqnarray}}
\newcommand{\beqno}{\begin{eqnarray*}}
\newcommand{\eeqno}{\end{eqnarray*}}

\def\b{{\bf b}}
\def\n{{\bf n}}

\begin{document}

\title*{On the sensitivity to model parameters in a filter stabilization technique for advection dominated advection-diffusion-reaction problems}
\titlerunning{A filter stabilization technique for ADR problems} 
\author{Kayla Bicol\inst{1} and Annalisa Quaini\inst{1}}
\authorrunning{Bicol and Quaini} 
\institute{
University of Houston, 4800 Calhoun Rd, Houston, TX 77204, USA, {\tt \{kmbicol,quaini\}@math.uh.edu}
This research has been supported in part by the NSF under grants DMS-1620384.
}
%
%
\maketitle

\abstract*{}


\abstract{We consider a filter stabilization technique with a deconvolution-based indicator function for the simulation of
advection dominated advection-diffusion-reaction (ADR) problems with under-refined meshes.
The proposed technique has been previously applied to the incompressible Navier-Stokes equations
and has been successfully validated against experimental data. However, it was found that 
some key parameters in this approach have a strong impact on the solution. 
To better understand the role of these parameters, 
we consider ADR problems, which are simpler than incompressible flow problems. 
For the implementation of the filter stabilization technique to ADR problems we adopt a three-step algorithm 
that requires (i) the solution of the given 
problem on an under-refined mesh,
(ii) the application of a filter to the computed solution, and (iii) a relaxation step. 
We compare our deconvolution-based approach to classical stabilization methods 
and test its sensitivity to model parameters on a 2D benchmark problem. 
}

\section{Introduction}\label{sec:intro}
We adapt to time-dependent advection-diffusion-reaction (ADR) problems
a filter stabilization technique proposed in \cite{Ervin_et_al2012} for evolution equations
and mostly developed for the Navier-Stokes equations \cite{Ervin_et_al2012,layton_CMAME,abigail_CMAME}.
This technique applied to the Navier-Stokes equations
has been extensively tested on both academic problems \cite{Ervin_et_al2012,layton_CMAME,abigail_CMAME}
and realistic applications \cite{BQV,Xu2018}. 
It was found in \cite{BQV} that key parameters in this approach have a strong impact on the solution. 
In order to understand the role of these key parameters, we 
apply the filter stabilization technique in the simplified context of ADR problems. 

It is well known that the Galerkin method for advection dominated ADR problems
can lead to unstable solutions with spurious oscillations \cite{B-quarteroniv2,C-hughes2,brooksh,hughesf4,douglasw1,hughes1,codina2}.
The proposed stabilization technique cures these oscillations by using an indicator function 
to tune the amount and location of artificial viscosity. The main advantage of this
 technique is that it can be easily implemented in legacy solvers.

\section{Problem definition}\label{sec:p_def}

We consider a time-dependent advection-diffusion-reaction problem defined on a bounded domain $\Omega \in \mathbb{R}^d$, 
with $d = 2, 3$, over a time interval of interest $(0, T]$:
\begin{align}
	&\partial_t u - \mu \Delta u  +  \nabla \cdot( \b u) + \sigma u = f \quad \text{in}~\Omega \times (0, T],  \label{eq:adr1}
\end{align}
endowed with boundary conditions:
\begin{align}
	&u = u_D \quad\quad\quad \text{on}~\partial \Omega_D \times (0, T], \label{eq:BC_D} \\
	& \mu \nabla u \cdot \n = g \quad \text{on}~\partial \Omega_N \times (0, T], \label{eq:BC_N}
\end{align}
and initial condition $u = u_0$ in $\Omega \times \{0\}$.
Here $\overline{\partial\Omega_D}\cup\overline{\partial\Omega_N}=\overline{\partial\Omega}$ and $\partial\Omega_D \cap\partial\Omega_N=\emptyset$ and $u_D$, $g$, and $u_0$ are given. 
In \eqref{eq:adr1}-\eqref{eq:BC_N}, $\mu$ is a diffusion coefficient, $\b$ is an advection field, 
$\sigma$ is a reaction coefficient, and $f$ is the forcing term.
For the sake of simplicity, we consider $\mu$, $\b$, and $\sigma$ constant.

Let $L$ be a characteristic macroscopic length for problem \eqref{eq:adr1}-\eqref{eq:BC_N}.
To characterize the solution of problem \eqref{eq:adr1}-\eqref{eq:BC_N}, we introduce
the P\'eclet number: $\displaystyle{\mathbb{P}\text{e} = || \b ||_\infty L/(2 \mu)}$.
We will assume that the problem is convection dominated, i.e. $|| \b ||_\infty \gg \mu$, which implies large P\'eclet numbers.
Notice that the role of the $\mathbb{P}\text{e}$ for advection-diffusion-reaction problems is similar to the 
role played by the Reynolds number for the Navier-Stokes equations.

In order to write the variational formulation of problem \eqref{eq:adr1}-\eqref{eq:BC_N}, 
we define the following spaces:
\begin{align}
	V &= \left\{ v: \Omega \rightarrow \R ,~v\in H^1({\Omega}),~v = u_D~\text{on}~\partial \Omega_D \right\}, \cl
	V_0 &= \left\{ v: \Omega \rightarrow \R ,~v\in H^1({\Omega}),~v = 0~\text{on}~\partial \Omega_D \right\}. \el
\end{align}
and bilinear form
\begin{align}
	b(u, w) =  (\mu \nabla u, \nabla w)_\Omega + (\nabla  \cdot (\b u), w)_\Omega + (\sigma u, w)_\Omega , \quad \forall u \in V, ~\forall w \in V_0. \label{eq:a}
\end{align}
The variational form of problem \eqref{eq:adr1}-\eqref{eq:BC_N} reads: Find $u \in V$ such that
\begin{align}\label{eq:weak}
	(\partial_t u,w)_\Omega + b(u, w) = (f, w)_\Omega + (g, w)_{\partial \Omega_N}, \quad \forall w \in V_0. 
\end{align}
The conditions for existence and unicity of the solution of problem \eqref{eq:weak} can be found, e.g., in 
\cite[Chapter~12]{B-quarteroniv2}. 

For the time discretization of problem \eqref{eq:weak}, we 
consider the backward Euler scheme for simplicity.
Let $\Delta t \in \mathbb{R}$, $t^n = n \Delta t$, with $n = 0, ..., N_T$ and $T = N_T \Delta t$. 
We denote by $y^n$ the approximation of a generic quantity $y$ at the time $t^n$. 
Problem \eqref{eq:weak} discretized in time reads: given $u^0 = u_0$, for $n \geq 0$ find
$u^{n+1} \in V$ such that
\begin{align}\label{eq:time_d}
	\frac{1}{\Delta t} (u^{n+1},w)_\Omega + b(u^{n+1}, w) =  \frac{1}{\Delta t} (u^{n},w)_\Omega + (f^{n+1}, w)_\Omega + (g^{n+1}, w)_{\partial \Omega_N},
\end{align}
$\forall w \in V_0$. We remark that time discretization approximates problem \eqref{eq:weak} by a sequence of quasi-static problems \eqref{eq:time_d}.
In fact, we can think of eq.~\eqref{eq:time_d} as the variational formulation of problem:
\begin{align}\label{eq:L}
	L_t u^{n+1} = {f}_t^{n+1},
\end{align}
where 
\begin{align}
L_t u^{n+1} = - \mu \Delta u^{n+1}  + \nabla \cdot (\b u^{n+1}) + \left( \sigma + \frac{1}{\Delta t} \right) u^{n+1}, ~{f}_t^{n+1} = \frac{u^{n}}{\Delta t} + f^{n+1}. \el
\end{align}

With regard to space discretization, we use the Finite Element method.
Let $\mathcal{T}_h = \{ K \}$ be a generic, regular Finite Element triangulation of the domain $\Omega$
composed by a set of finite elements, indicated by $K$. 
As usual, $h$ refers to the largest diameter of the elements of $\mathcal{T}_h$.
Let $V_h$ and $V_{0,h}$ be the finite element spaces approximating $V$ and $V_0$, respectively.
The fully discrete problem reads:
given $u^0_h$, for $n \geq 0$ find $u^{n+1}_h \in V_h$ such that
\begin{align}\label{eq:spacetime_d}
	\frac{1}{\Delta t} (u^{n+1}_h,w)_\Omega + b(u^{n+1}_h, w) =  \frac{1}{\Delta t} (u^{n}_h,w)_\Omega + (f^{n+1}_h, w)_\Omega + (g^{n+1}_h, w)_{\partial \Omega_N}, 
\end{align}
$\forall w \in V_{0,h}$, where $u^0_h$, $f^{n+1}_h$, and $g^{n+1}_h$ are appropriate finite element approximations
of $u_0$, $f^{n+1}$, and $g^{n+1}$, respectively. It is known that the 
solution of problem \eqref{eq:spacetime_d} converges optimally to the solution
of problem \eqref{eq:weak}. However, 
the Finite Element method can perform poorly if the coerciveness
constant of the bilinear form \eqref{eq:a} is small in comparison with its continuity constant.
In particular, the error estimate can have a very large multiplicative
constant if $\mu$ is small with respect to $|| \b ||_\infty$, i.e. when $\mathbb{P}\text{e}$
is large. In those cases, 
the finite element solution $u_h$ can be globally polluted with strong spurious oscillations.
To characterize the solution of problem \eqref{eq:spacetime_d}, we introduce
the local counterpart of the P\'eclet number: $\displaystyle{\mathbb{P}\text{e}_h = || \b ||_\infty h/(2 \mu)}$.

Several stabilization techniques have been proposed to eliminate, or at least reduce, the 
numerical oscillations produced by the standard Galerkin method in case of large $\mathbb{P}\text{e}$. 
In the next section, we will go over a short review of these stabilization techniques before
introducing our filter stabilization method in Sec.~\ref{eq:filter_stab}.

\subsection{Overview of stabilization techniques}\label{sec:overview}

We will restrict our attention to stabilization techniques that consists of 
adding a stabilization term $b_s(u^{n+1}_h, w)$ to the left-hand side of time-discrete problem \eqref{eq:time_d}.
In the following, we will use $\tau$ to denote a stabilization parameter that can depend on the element
size $h$ and the equation coefficients. Parameter $\tau$ takes different values for the different stabilization schemes.  
We will use the broken inner product $(\cdot, \cdot)_K = \sum_{K} ( \cdot, \cdot)$, where $\sum_{K}$
denotes summation over all the finite elements.

Perhaps the easiest way to stabilize problem \eqref{eq:time_d} 
is by introducing  artificial viscosity either in the whole domain,
leading to $b_s(u^{n+1}_h, w)= (\tau \nabla u_h^{n+1}, \nabla w)_K$,
or streamwise, leading to $b_s(u^{n+1}_h, w) = (\tau \b \cdot \nabla u_h^{n+1}, \b \cdot \nabla w)_K$.
See, e.g., \cite{C-hughes2, B-quarteroniv2}. 
In this way, the effective $\mathbb{P}\text{e}_h$ becomes smaller. 
The artificial viscosity $\tau$ in these schemes is proportional to $h$. 
The drawbacks of these schemes are that they are first order accurate only
and not strongly consistent.
An improvement over the artificial viscosity schemes is given by the strongly consistent 
stabilization methods. In fact, strong consistency allows the stabilized method to maintain
the optimal accuracy. 

Let us introduce the residual for problem \eqref{eq:L} and the skew-symmetric part of operator $L_t$:
\begin{align}
	R(u_h^{n+1})  = {f}_t^{n+1} - L_t u_h^{n+1}, \quad L_{SS} v  = \frac{1}{2} \div (\b v) + \frac{1}{2} \b \cdot \nabla v. \el
\end{align}
One of the most popular strongly consistent stabilized finite element methods is the 
Streamline Upwind Petrov-Galerkin (SUPG) method \cite{brooksh}, for which 
$ b_s(u^{n+1}_h, w) = - (\tau R(u_h^{n+1}), L_{SS} w)_K$.
The Galerkin Least Squares (GLS) method \cite{hughesf4}
is a generalization of the SUPG method: $ b_s(u^{n+1}_h, w) = - (\tau R(u_h^{n+1}), L_t w)_K$.
The Douglas-Wang method \cite{douglasw1} replaces
$ L_t w$ in the GLS method with $-L_t^*w$, where $L_t^*$ is the adjoint of operator $L_t$. Thus, we have:
$ b_s(u^{n+1}_h, w) = (\tau R(u_h^{n+1}), L_t^*w)_K$. For the SUPG, GLS, and Douglas-Wang
methods, $\tau = \delta h_K/| \b |$ for $\delta > 0$. 
Finally, we mention a method based on the Variational Multilscale approach \cite{hughes1}
called algebraic subgrid scale (ASGS), for which $b_s(u^{n+1}_h, w) = (\tau R(u_h^{n+1}), L_t^*w)_K$.
The difference between the Douglas-Wang and ASGS method consists in the choice for parameter $\tau$.
One possibility for $\tau$ in the ASGS method is $\tau = [ 4 \mu/h_K^2 + 2 |\b|/h_K + \sigma ]^{-1}$ (see
\cite{codina2}).

We report in Table \ref{tab:stab} a summary of the methods in this overview. All of the strongly consistent
stabilization techniques come with stability estimates that improve the one that can be obtained for the Galerkin method.
See \cite{brooksh,hughesf4,douglasw1,B-quarteroniv2,hughes1,codina2}.

\begin{table}
\caption{Stabilization term $b_s(u^{n+1}_h, w)$ for some stabilization methods.}
\label{tab:stab} 
	\begin{center}
		\begin{tabular}{c|c|c}
			\hline
			&Stabilization method &  $b_s(u^{n+1}_h, w)$  \\
			\hline
			\multirow{2}{6em}{Not strongly consistent} & Artificial viscosity & $(\tau \nabla u_h^{n+1}, \nabla w)_K$ \\
			& Streamline upwind & $(\tau \b \cdot \nabla u_h^{n+1}, \b \cdot \nabla w)_K$ \\
			\hline
			\multirow{4}{6em}{Strongly consistent} & Streamline upwind Petrov-Galerkin & $- (\tau R(u_h^{n+1}), L_{SS} w)_K$ \\
			& Galerkin least-squares & $-(\tau R(u_h^{n+1}), L_t w)_K$ \\
			& Douglas-Wang  & $(\tau R(u_h^{n+1}), L_t^*w)_K$ \\
			& Variational Multiscale, ASGS  & $(\tau R(u_h^{n+1}), L_t^*w)_K$ \\
			\hline
		\end{tabular}
	\end{center}
\end{table}


\section{A filter stabilization technique}\label{eq:filter_stab}

We adapt to the time-dependent advection-diffusion-reaction problem defined in Sec.~\ref{sec:p_def}
a filter stabilization technique proposed in \cite{Ervin_et_al2012}.
For the implementation of this stabilization technique we adopt an algorithm called
\emph{evolve-filter-relax} (EFR) that was first presented in \cite{layton_CMAME}. 
The EFR algorithm applied to problem \eqref{eq:L} with boundary conditions
\eqref{eq:BC_D}-\eqref{eq:BC_N} reads: given $u^n$

\begin{enumerate}[label=(\roman*)]
	\item \emph{Evolve}: find the intermediate solution $v^{n+1}$ such that
	\begin{align}
		&L_t v^{n+1} = f_t^{n+1} \quad \text{in}~ \Omega  , \label{eq:evolve} \\
		&v^{n+1} = u_D \quad\quad\quad \text{on}~\partial \Omega_D , \label{eq:evolve-2}\\
		& \mu \nabla v^{n+1} \cdot \n = g \quad \text{on}~\partial \Omega_N. \label{eq:evolve-3}
	\end{align}
	\item \emph{Filter}: find $\overline{v}^{n+1}$ such that
	\begin{align}
		&\overline{v}^{n+1} - \delta^2 \nabla \cdot ( a(v^{n+1}) \nabla \overline{v}^{n+1}) = v^{n+1} &&\text{in}~ \Omega  , \label{eq:filter}
		\\
		&\overline{v}^{n+1} = u_D &&\text{on}~\partial \Omega_D , \label{eq:filter-2} \\
		& \mu \nabla \overline{v}^{n+1} \cdot \n = 0 &&\text{on}~\partial \Omega_N.  \label{eq:filter-3}
	\end{align}
	Here, $\delta$ can be interpreted as the \emph{filtering radius} 
	(that is, the radius of the neighborhood were the filter extracts information)
	and $a(\cdot) \in (0,1]$ is a scalar function called \emph{indicator function}. The indicator function has to be
	such that $a(v^{n+1})\simeq 1$ where $v^{n+1}$ does need to
	be filtered from spurious oscillations, and
	$a(v^{n+1})\simeq 0$ where $v^{n+1}$ does not need to
	be filtered.
	\item \emph{Relax}: set 
	\begin{align}
		u^{n+1} = (1 - \chi) v^{n+1} + \chi \overline{v}^{n+1}, \label{eq:relax}
	\end{align}
	where $\chi\in (0,1]$ is a relaxation parameter.
\end{enumerate}


The EFR algorithm has the advantage of modularity: since
the problems at steps (i) and (ii) are numerically standard, they can be solved with legacy solvers
without a considerable implementation effort.
Algorithm \eqref{eq:evolve}-\eqref{eq:relax} is sensitive to the choice of key parameters
$\delta$ and $\chi$ \cite{BQV,BQRV}.
A common choice for $\delta$ is $\delta = h$.
However, 
in \cite{BQV} it is suggested that taking $\delta = h$ might lead to to excessive numerical diffusion
and it is proposed to set $\delta = h_{min}$, where $h_{min}$ is the length of the shortest edge in the mesh.
As for $\chi$, in \cite{layton_CMAME} the authors  support the choice $\chi=O(\Delta t)$
because it guarantees that the numerical dissipation vanishes as $h \rightarrow 0$ regardless of $\Delta t$.
In \cite{BQV} the value of $\chi$ is set with a heuristic formula which depends on both physics
ans discretization parameters.

Different choices of $a(\cdot)$ for the Navier-Stokes equations have been 
proposed and compared in  \cite{layton_CMAME,O-hunt1988,Vreman2004,Bowers2012}.
Here, we focus on a class of deconvolution-based indicator functions:
\begin{equation}
a(u) = a_{D}(u) = \left|  u - D (F(u)) \right|, \label{eq:a_deconv}
\end{equation}
where $F$ is a linear, invertible, self-adjoint, compact operator from a Hilbert space $V$ to itself,
and $D$ is a bounded regularized approximation of $F^{-1}$. In fact, since 
 $F$ is compact, the inverse operator $F^{-1}$ is unbounded.
The composition of the two operators $F$ and $D$ can be interpreted as a low-pass filter.

A possible choice for $D$ is the Van Cittert deconvolution operator $D_N$, defined as
\begin{equation}
D_N = \sum_{n = 0}^N (I - F)^n. \el
\end{equation}
The evaluation of $a_D$ with $D=D_N$ (deconvolution of order $N$) requires then to apply the filter $F$ a total of $N+1$ times. 
Since $F^{-1}$ is not bounded, in practice $N$ is chosen to be small, as the result of a trade-off between accuracy (for a regular solution) and filtering (for a non-regular one).

We select $F$ to be the linear Helmholtz filter operator $F_H$ \cite{germano} defined by
\begin{equation}
F=F_H \equiv \left(I - \delta^2 \Delta \right)^{-1}. \el
\end{equation}
It is possible to prove  \cite{Dunca2005} that
\begin{align}
	u - D_N (F_H(u)) = (-1)^{N+1} \delta^{2N+2} \Delta^{N+1} F_H^{N+1} u. \label{eq:duncan_eps}
\end{align}
Therefore, $a_{D_N}(u)$ is close to zero in the regions of the domain where $u$ is smooth. Indicator function (\ref{eq:a_deconv}) with $D=D_N$ and $F=F_H$ has been proposed in \cite{abigail_CMAME} for the Navier-Stokes equations.
Algorithm \eqref{eq:evolve}-\eqref{eq:relax} with indicator function (\ref{eq:a_deconv}) 
is also sensitive to the choice of $N$ \cite{BQV,BQRV}.

In order to compare our approach with the stabilization techniques 
reported in Sec.~\ref{sec:overview}, let us assume that problem \eqref{eq:adr1} is supplemented 
with homogeneous Dirichlet
boundary conditions on the entire boundary, i.e. $\partial \Omega_D = \partial \Omega$ and
$u_D = 0$ in \eqref{eq:BC_D}. Let us start by writing the weak form of eq.~\eqref{eq:evolve}:
\begin{align}
	(L_t v^{n+1},w) = \left(f_t^{n+1},w \right)_\Omega. \label{eq:evolve_weak}
\end{align}
Next, we apply operator $L_t$ to eq.~\eqref{eq:filter}
and write the corresponding weak form, using also eq.~\eqref{eq:evolve_weak}:
\begin{align}
	(L_t \overline{v}^{n+1},w) - ( \nabla \cdot ( \overline{\mu} \nabla \overline{v}^{n+1}),L_t^*w) = \left(f_t^{n+1},w \right)_\Omega, \quad \overline{\mu} = \delta^2 a(v^{n+1}). \label{eq:Lfilter_weak}
\end{align}
Here, $\overline{\mu}$ is the artificial viscosity introduced by our stabilization method. 
Now, we multiply eq.~\eqref{eq:evolve_weak} by $(1 -\chi)$ and add it to eq.~\eqref{eq:Lfilter_weak}
multiplied by $\chi$. Using the relaxation step \eqref{eq:relax}, we obtain:
\begin{align}
	(L_t u^{n+1},w) - \chi (\nabla \cdot ( \overline{\mu}  \nabla \overline{v}^{n+1})),L_t^*w) = \left(f_t^{n+1},w \right)_\Omega, \label{eq:e_f_r}
\end{align}

The second term at the left-hand side in \eqref{eq:e_f_r} is the stabilization term added by the filter stabilization
technique under consideration.
Notice that eq.~\eqref{eq:e_f_r} is a consistent perturbation of the original advection-diffusion-reaction problem.
The perturbation vanishes with coefficient $\chi$, which goes to zero with the discretization parameters 
(recall that a possible choice is $\chi=O(\Delta t)$).
Using \eqref{eq:relax} once more, we can rewrite the stabilization term as:
\begin{align}
	b_s(u^{n+1},w) = -  ( \nabla \cdot ( \overline{\mu} \nabla u^{n+1}),L_t^*w)  + (1 - \chi) ( \nabla \cdot ( \overline{\mu}\nabla v^{n+1}),L_t^*w).  
\end{align}
We see that,
the stabilization term here does not depend only on the end-of-step solution $u^{n+1}$.
We remind that usually $\delta = h$.
Thus, as $h \rightarrow 0$ the artificial viscosity $\overline{\mu}$ in \eqref{eq:Lfilter_weak}
vanishes. It is then easy to see that the filter stabilization technique we consider
is consistent, although not strongly. 

\section{Numerical results}

We consider a benchmark test proposed in \cite{Volker_Schmeyer2008}. 
The prescribed solution is given by:
\begin{align}
u(x,y,t) = &16 \sin(\pi t) x(1-x)y(1-y) \cl
&\cdot \left[ \frac{1}{2} + \frac{\arctan(2 \mu^{-1/2} (0.25^2 - (x - 0.5)^2 - (y - 0.5^2)))}{\pi}\right], \label{eq:exact_sol}
\end{align}
in $\Omega = (0,1) \times (0,1)$ and in time interval $(0, 0.5]$.
We set $\sigma = 1$,
$\mu = 10^{-5}$, and $\b = [2,3]^T$,  which yield 
$\mathbb{P}\text{e} = 1.5\cdot 10^{5}$. 
Solution \eqref{eq:exact_sol} is a hump changing its height in the course of the time. 
The internal layer in solution \eqref{eq:exact_sol} has size $O(\sqrt{\mu})$.
The forcing term $f$ in \eqref{eq:adr1} and the initial condition $u_0$ follow from \eqref{eq:exact_sol}. 
We impose boundary condition \eqref{eq:BC_D} with $u_D = 0$ on the entire boundary, which is 
consistent with exact solution \eqref{eq:exact_sol}. We use this test to compare the solution
computed by the EFR method with the solution given by other methods, and to show
the sensitivity of the solution computed by the EFR method to parameters $N$ and $\delta$. 
The sensitivity to $\chi$ will be object of future work.  
All the computational results have been obtained
with FEniCS \cite{fenics,LoggMardalEtAl2012a,AlnaesBlechta2015a}.

We take $\Delta t = 10^{-3}$. 
We consider structured meshes with 5 different refinement levels $\ell = 0, \cdots, 4$
and $\mathbb{P}_2$ finite elements.
Triangulation $\mathcal{T}_{h_\ell}$ of $\Omega$ consists of $n_\ell^2$ sub-squares, each of which
is further divided into 2 triangles. The associated mesh size is $h_\ell = \sqrt{2}/n_\ell$.
In Table \ref{tab:meshes}, we report $n_\ell$ and $\mathbb{P}\text{e}_h$
for each of the meshes under consideration. We see that even on the finest mesh, 
the local P\'eclet number is much larger than 1. Table \ref{tab:meshes} gives also
the value of $\chi$ used for the EFR method on the different meshes. 

\begin{table}
\caption{Number of partitions $n_\ell$ for each side, local P\'eclet number $\mathbb{P}\text{e}_h$,
and value of $\chi$ for the EFR method for the meshes associated to 5 different refinement levels.}\label{tab:meshes} 
	\begin{center}
		\begin{tabular}{c|c|c|c|c|c}
			\hline
			 Refinement level  & $\ell = 0$ & $\ell = 1$ & $\ell = 2$ &  $\ell =  3$ & $\ell =  4$\\
			\hline
			$n_\ell$ & 25 & 50 & 100 & 200 & 400  \\
			\hline
			$\mathbb{P}\text{e}_h$ & 8485.3 & 4242.6 & 2121.3 & 1060.7 & 530.3  \\
			\hline
			$\chi$ & 1 & 1/2 & 1/4 & 1/16 & 1/256  \\
			\hline
		\end{tabular}
	\end{center}
\end{table}

\begin{figure}
\begin{center}
\includegraphics[height=.4\textwidth]{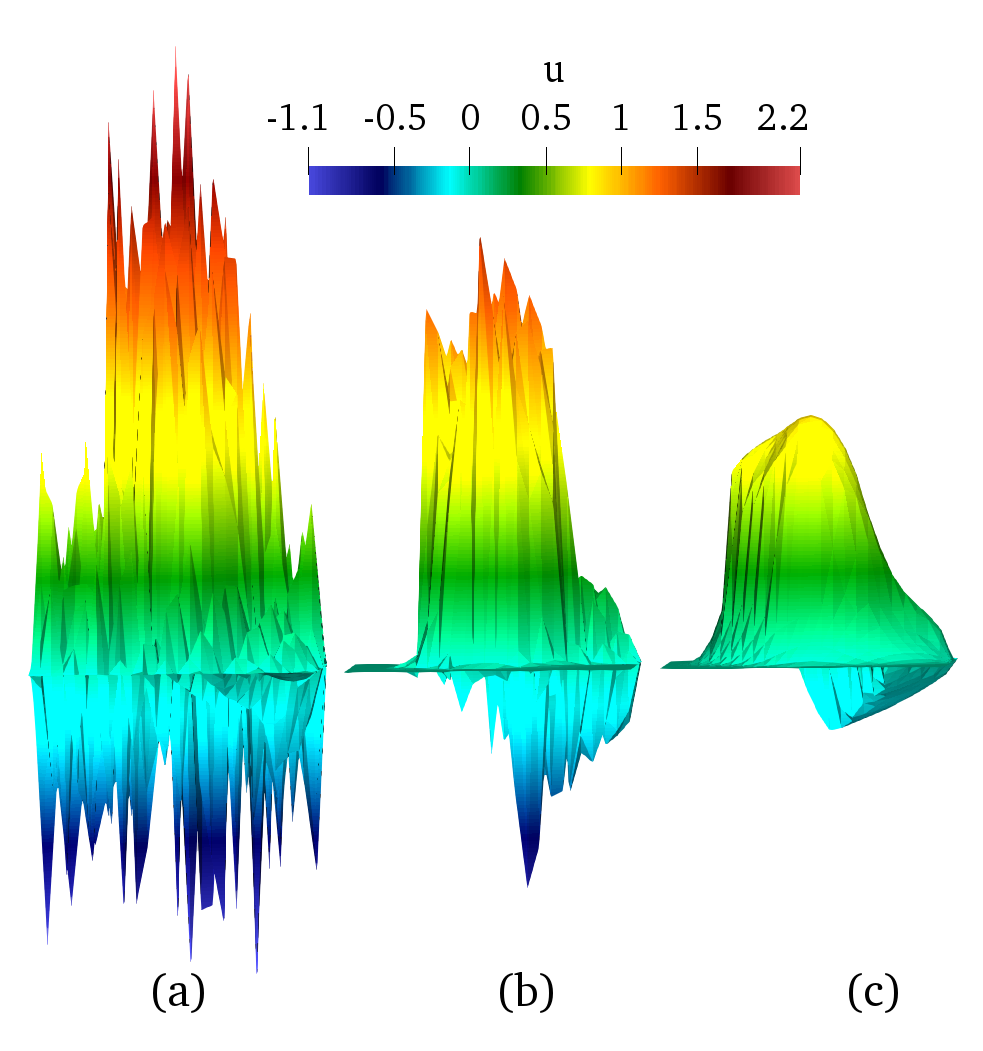}~
\includegraphics[height=.35\textwidth]{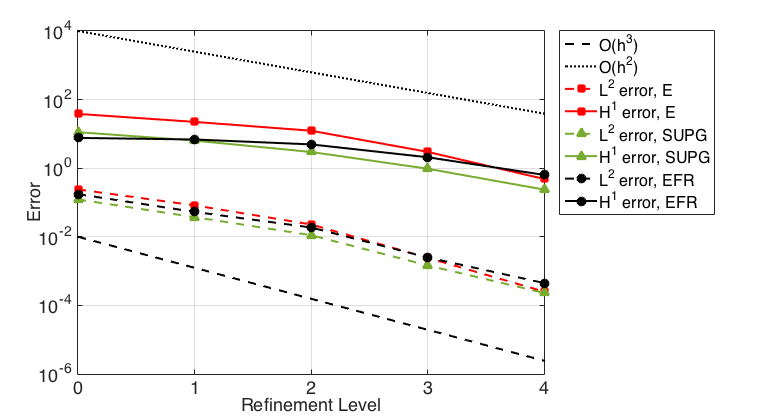}
\end{center}

\caption{Left:
Solution at $t = 0.5$ computed on mesh $\ell = 0$ with (a) the standard Galerkin method, 
(b) the SUPG method, and (c) the EFR method with $\delta = 1/n_\ell$ and $N =0$. 
Right: $L^2$ and $H^1$ norms of the error for $u$ at $t = 0.5$
given by the standard Galerkin method (E), the SUPG method, and the EFR method
plotted against the refinement level.
}
\label{fig:comparison}       
\end{figure}

Since the problem is convection-dominated and
the solution has a (internal) layer, the use of a stabilization method is necessary. 
See Fig.~\ref{fig:comparison} (left) for a comparison of the solution at $t = 0.5$ computed on mesh $\ell = 0$
with the standard Galerkin element method, the SUPG method, and 
the EFR method with $\delta = 1/n_\ell$ and $N =0$.
We see that the solution obtained with the non-stabilized 
method is globally polluted with spurious oscillations. Oscillations are still present in 
SUPG method (mainly at the top of the hump and in the right upper part of the domain), 
but they are reduced in amplitude. The amplitude of the oscillations is further reduced
in the solution computed with the EFR method. The other strongly consistent stabilization methods reported in
Sec.~\ref{sec:overview} give results very similar to the SUPG methods. For this reason, 
those results are omitted. Fig.~\ref{fig:comparison} (right) shows the $L^2$ and $H^1$ norms
of the error for the solution at $t = 0.5$ plotted against the mesh refinement level $\ell$.
We observe that the EFR method gives errors comparable to those given by the SUPG method
on the coarser meshes.
When the standard Galerkin method gives
a smooth approximation of the solution, e.g. with mesh $\ell = 4$, 
the errors given by the EFR method are comparable to the errors given by the
standard Galerkin method.

\begin{figure}
\begin{center}
\includegraphics[height=.35\textwidth]{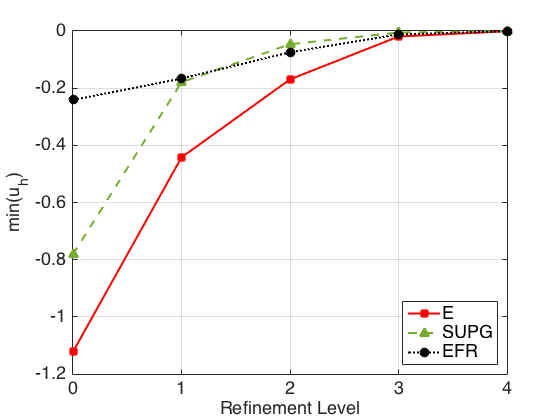}~
\includegraphics[height=.35\textwidth]{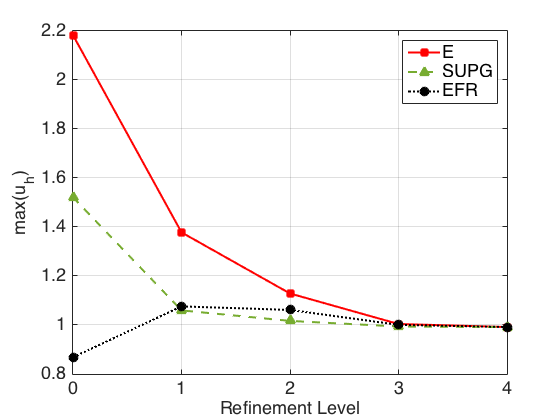}
\end{center}
\caption{Minimum value (left) and maximum value (right) of the 
solution at $t = 0.5$ computed by the standard Galerkin method (E), 
the SUPG method, and the EFR method with $\delta = 1/n_\ell$, $N =0$, and $\chi = 1$
plotted against the refinement level.
}
\label{fig:minmax}       
\end{figure}

In Fig.~\ref{fig:minmax}, we report the minimum value (left) and maximum value (right) of the 
solution at $t = 0.5$ computed by the standard Galerkin method, the SUPG method, 
and EFR method with $\delta = 1/n_\ell$ and $N =0$. 
plotted against the refinement level $\ell$. 
The SUPG method does not eliminate the under- and over-shoots given by
a standard Galerkin method but it reduces their amplitude. It is well-known that all finite element methods
that rely on streamline diffusion stabilization produce under- and over-shoots 
in regions where the solution gradients are steep and not aligned with the direction of $\b$. 
From Fig.~\ref{fig:minmax}, we see that 
the EFR method gives under- and overs-shoots 
of smaller or comparable amplitude when compared to the SUPG method. 
In some practical applications, such imperfections are small in magnitude and can be tolerated. 
In other cases, it is essential to ensure that
the numerical solution remains nonnegative and/or devoid of spurious oscillations.
This can be achieved with, e.g., discontinuity-capturing or shock-capturing techniques
\cite{HUGHES1986329,CODINA1993325,DESAMPAIO20016291,Kuzmin2010}. 
However, this is outside the scope of the present work. 

\begin{figure}
\begin{center}
\includegraphics[height=.31\textwidth]{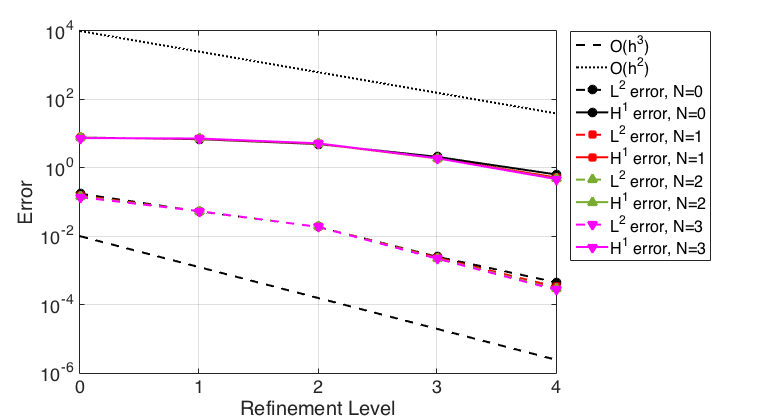}~
\includegraphics[height=.31\textwidth]{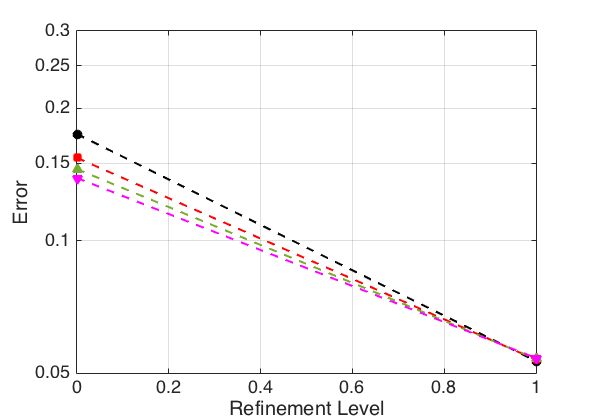}
\end{center}
\caption{
Left: $L^2$ and $H^1$ norms of the error for $u$ at $t = 0.5$
given by the EFR method with $\delta = 1/n_\ell$ and
$N = 0, 1, 2, 3$ plotted against the refinement level.
Right: zoomed-in view around $\ell = 0, 1$.}
\label{fig:different_N}       
\end{figure}

Next, we focus on the EFR algorithm and vary the order of the deconvolution $N$. 
In Fig.~\ref{fig:different_N} (left) we show
$L^2$ and $H^1$ norms of the error for $u$ at $t = 0.5$
given by the EFR method with $\delta = 1/n_\ell$ and
$N = 0, 1, 2, 3$ plotted against the refinement level.
The only visible difference when $N$ varies is for the finer meshes, with 
both errors slightly decreasing as $N$ is increased.
Fig.~\ref{fig:different_N} (right) displays a zoomed-in view of Fig.~\ref{fig:different_N} (left)
around $\ell = 0, 1$.
It shows that also for the coarser meshes 
the errors get slightly smaller when $N$ increases. We recall that 
indicator function \eqref{eq:a_deconv} with $D = D_N$ requires to apply the 
Helmholz filter $N+1$ times. So, the slightly smaller errors for large $N$ 
come with an increased computational time. 

\begin{figure}
\begin{center}
\includegraphics[width=.7\textwidth]{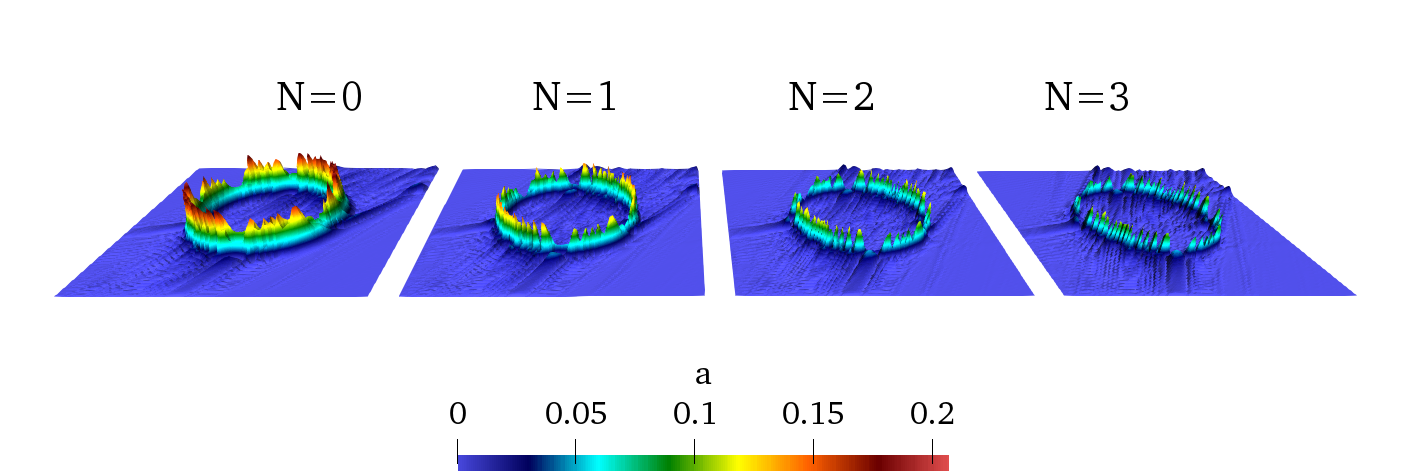}
\end{center}
\caption{Indicator function on mesh $\ell = 2$ at $t = 0.5$ for $\delta = 1/n_\ell$ and
$N = 0, 1, 2, 3$.}
\label{fig:ind_nx100}    
\end{figure}

In Fig.~\ref{fig:ind_nx100}, we report the indicator function at $t = 0.5$ for $\delta = 1/n_\ell$ and
$N = 0, 1, 2, 3$. In all the cases, the largest values of the indicator
function are around the edge of the hump. Moderate values are aligned
with the direction of $\b$. While we see some differences in the indicator functions for $N = 0$ and $N=1$, 
for $N > 1$ the indicator function does not seem to change substantially.

\begin{figure}
\begin{center}
\includegraphics[height=.31\textwidth]{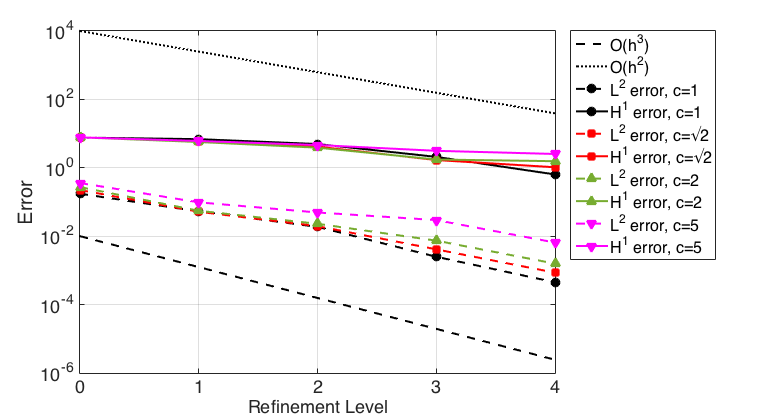}~
\includegraphics[height=.33\textwidth]{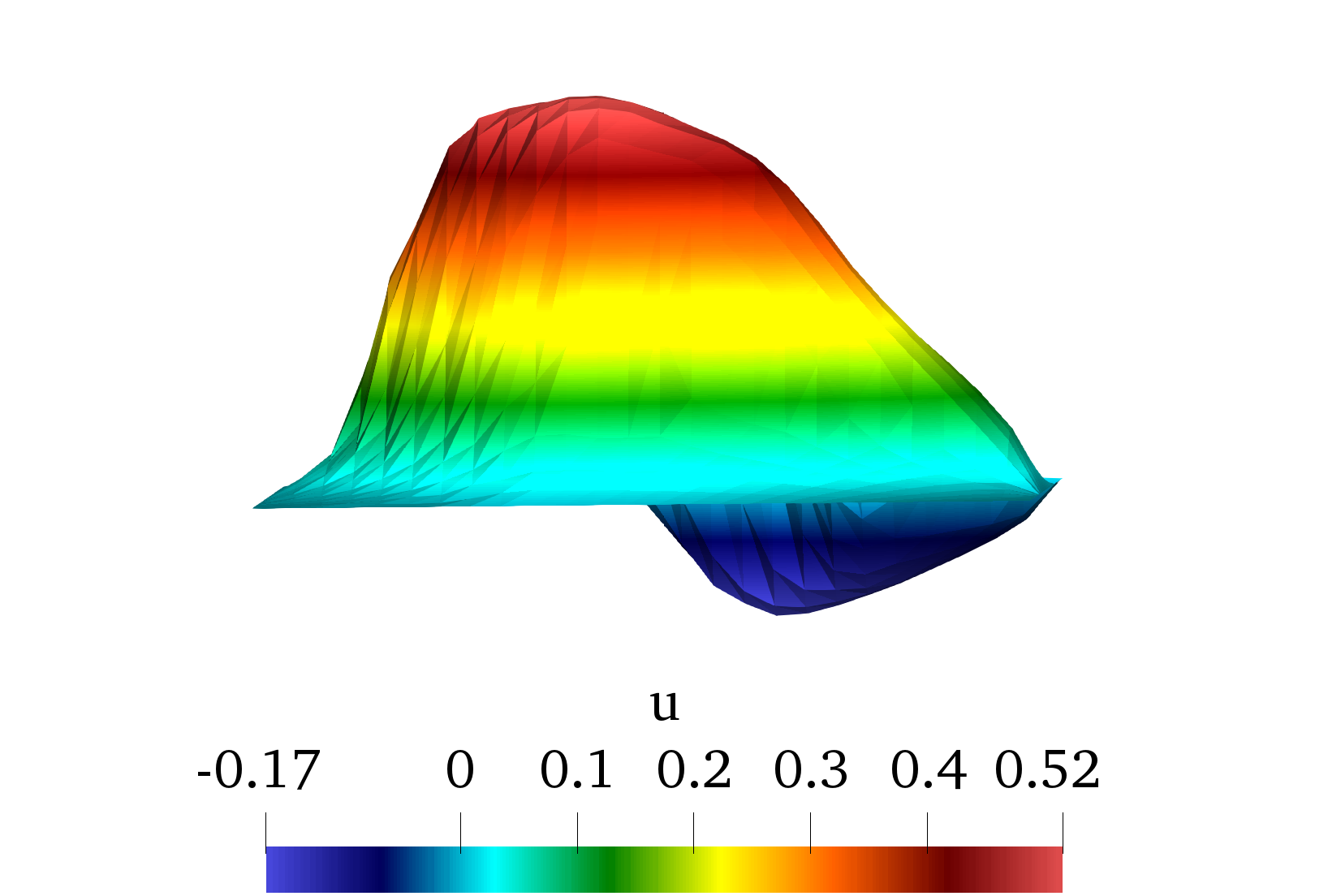}
\end{center}
\caption{Left:
$L^2$ and $H^1$ norms of the error for $u$ at $t = 0.5$
given by the EFR method with $N =0$ and $\delta = c/n_\ell$,
$c = 1, \sqrt{2}, 2, 5$, plotted against the refinement level.
Right: solution at $t = 0.5$ computed by the EFR method with 
$\delta = 5/n_\ell$ on mesh $\ell = 0$.
}
\label{fig:different_delta}       
\end{figure}

Finally, we fix $N = 0$ and vary $\delta$.
In Fig.~\ref{fig:different_delta} (left) we show
$L^2$ and $H^1$ norms of the error for $u$ at $t = 0.5$
given by the EFR method with
$\delta = c/n_\ell$, $c = 1, \sqrt{2}, 2, 5$, plotted against the refinement level. 
Notice that $c = \sqrt{2}$ corresponds to the choice $\delta = h$, 
while $c = 1$ corresponds to $\delta = h_{min}$.
In \cite{BQV}, it was found that $\delta = h_{min}$ makes the numerical results 
(for a Navier-Stokes problem on unstructured meshes) in better agreement 
with experimental data. Our results confirm that $\delta = h_{min}$ is
the best choice. In fact, it minimizes the error and gives optimal convergence rates. 
Higher values of $c$, i.e. $c > 1$, seem to spoil the convergence rate of the EFR method. 
From Fig.~\ref{fig:different_N} (left) and \ref{fig:different_delta} (left) we see 
that the computed solution is much more sensitive to $\delta$ than it is to $N$. 
Fig.~\ref{fig:different_delta} (right) shows the solution at time $t = 0.5$ computed  
by the EFR method with $\delta = 5/n_\ell$ on mesh $\ell = 0$.
Remember that $\delta$ is the filtering radius, i.e.~the radius of the circle over 
which we average (in some sense) the 
solution. Thus, it is not surprising that for a large value of $\delta$ the EFR method
has an over-smoothing effect. 

\section{Conclusions}

We considered a deconvolution-based filter stabilization technique
recently proposed for the Navier-Stokes equations and adapted it to the numerical solution of
advection dominated advection-diffusion-reaction problems with under-refined meshes.
Our stabilization technique is consistent, although not strongly. 
For the implementation of our approach we adopted a three-step algorithm 
called evolve-filter-relax (EFR) that can be easily realized within a legacy solver.
We showed that the EFR algorithm is competitive when compared to classical 
stabilization methods on a benchmark problem that features an analytical solution.
However, special care has to be taken in setting the filtering radius $\delta$ 
in order to avoid over-smoothing.

\bibliographystyle{plain}
\bibliography{latexbi,references,lifev} 

\end{document}